\documentstyle[12pt]{amsart}
\newtheorem{theorem}{Theorem}[section]

\newtheorem{lemma}[theorem]{Lemma}
\newtheorem{proposition}[theorem]{Proposition}

\numberwithin{equation}{section}

\def\be{\begin{equation}}
\def\en{\end{equation}}
\def\ba{\begin{eqnarray}}
\def\ea{\end{eqnarray}}

\def\e1{\epsilon}
\def\o1{\omega}
\def\01{\Omega}
\def\al{\alpha}
\def\ga{\gamma}
\def\g1{\Sigma}
\def\lmd{\lambda}
\def\l1{\Lambda}
\def\v1{\varphi}
\def\d1{\delta}
\def\f1{\frac}
\def\h2{{\bf H}}
\def\a2{{\bf A}}
\def\x2{{\bf X}}
\def\t1{\theta}
\def\b1{\beta}

\def\bs{\begin{eqnarray*}}
\def\es{\end{eqnarray*}}
\def\p{\partial}
\def\m1{\Theta}
\def\w1{\wedge}

\def\ctg{\cot}
\def\tg{\tan}

\begin{document}

\title[Harmonic Hopf Constructions]{Harmonic Hopf Constructions between spheres II}

\author{Weiyue Ding, Huijun Fan, Jiayu Li}

\thanks{The research was supported in part by  NSF of China,
National Key Basic Research Fund (1999075109), and by the
Outstanding Young Scientists Grants (199925101). The second author also thanks MPI in Leipzig for partially supporting.}

\address{Department of Mathematics\\Beijing University
\\Beijing 100080\\P. R. of China}
\email{dingwy@@math.pku.edu.cn}

\address{Institute of Mathematics\\
 Academia Sinica\\ Beijing 100080\\P. R. of
China} \email{fanhj@@Math08.math.ac.cn}

\address{Institute of Mathematics\\
Fudan University, Shanghai and
Academia Sinica, Beijing\\P.R. China}
\email{lijia@@math03.math.ac.cn}

\date{}



\maketitle

\section{Introduction and the main result}
This paper can be seen as the final remark of the previous paper
written by the first author [D]. We consider the existence of
harmonic maps between two spheres, via Hopf constructions. Given a
non trivial bi-eigenmap f : $S^p \times S^q \longrightarrow S^n$
with bi-eigenvalue ($\lambda ,\mu$) ($\lambda ,\mu >0$) and a
continuous function $\alpha$ : [ 0, $\frac{\displaystyle
\pi}{\displaystyle 2}$ ]$\longrightarrow$ [ 0, $\displaystyle \pi$
] with $\alpha(0)=0$ , $\alpha(\frac {\displaystyle
\pi}{\displaystyle 2})$ = $\displaystyle \pi$ , one defines a map
u: $S^{p+q+1} \longrightarrow S^{n+1}$ , called the $\alpha$-Hopf
construction on f, by $$u ( \sin t\cdot x , \cos t\cdot y )=(
\sin\alpha(t)f(x,y) , \cos \alpha (t) ),$$ where $x \in S^p , y
\in
 S^q$ and $t \in [ 0 , \frac{\displaystyle\pi}{\displaystyle 2} ].$
 It is known [ER] that u is a harmonic map if and only if $\alpha$
 is a solution of the o.d.e.
 \begin{equation}\label{(1)}\ddot{\alpha}+(p\cdot \ctg t - q\cdot \tg t)\dot{\alpha}-
 (\frac{\lambda}{\sin^2t}+\frac{\mu}{\cos^2t})
 \sin\alpha\cdot \cos\alpha=0\end{equation}with the boundary value
 condition
 \begin{equation}\label{(2)}\lim_{t\longrightarrow0^+}\alpha(t)=0,
 \lim_{t\longrightarrow{\frac{\pi}{2}-0}}\alpha(t)
 =\pi\end{equation}
Ding [D] proves that , if p22 and q22 , then
(\ref{(1)})-(\ref{(2)}) has a solution .In the case that p=q=1 ,
Eells-Ratto [ER] prove that a (1.1)-(1.2) is solvable if and only
if $\lambda=\mu$. In  this paper ,we consider the remaining case ,
i.e. p=1 and $q\geq 2$. In this case , it is proved in [ER] that ,
a necessary condition for (1.1)-(1.2) to be solvable is that
$q\lambda< \mu$. Ding([D] Remark 1.2) conjectured that, it is also
a sufficient condition. In this paper, we prove the conjecture.

 \begin{theorem}  If p=1, $q>1, \lambda\geq1$  and $\mu > \lambda
 q$, then the prob.(1.1)-(1.2)
 has a solution $\alpha$  with  $0<\alpha(t)<\pi$  for
 $t\in(0,\frac{\pi}{2})$.
\end{theorem}

\noindent
{\bf Some applications}\\
\indent
An immediate consequence of Thoerem 1.1 is that it can provide non-contractible harmonic maps between spheres. For example, $\phi=e^{i\lambda\theta}$ is an eigenmap with eigenvalue $\lambda^2$ and if $\psi:\;S^q\to S^{2m-1}$ is another eigenmap with eigenvalue $\mu$, then the map $F=\phi\psi:\;S^1\times S^q\to S^{2m-1}\subset {\Bbb C}^m$ is a bi-eigenmap with bi-eigenvalue $(\lambda^2,\mu)$. If $\mu>\lambda^2 q$, then there exists a harmonic map from $S^{q+2}$ to $S^{2m}$.\\ 
\indent
Using orthogonal multiplication and Hopf constructions in the following way, we can obtain more information about the homotopy groups of spheres, which is based on Theorem 1.1.\
\indent
An orthogonal multiplication is a bilinear map $f:\;{\Bbb R}^k\times{\Bbb R}^l\to{\Bbb R}^n$ such that 
$$
|f(x,y)|=|x||y|,\;\forall x\in {\Bbb R}^k,\;y\in{\Bbb R}^l.
$$
A Hopf construction on $f$ is a map
$$
F_f:\;{\Bbb R}^k\times{\Bbb R}^l\rightarrow {\Bbb R}^{n+1}
$$
which sends $(x,y)$ to $(2f(x,y),|x|^2-|y|^2)$.\\
Its restriction induces a map
$$
\psi:\;S^{k+l-1}\rightarrow S^n
$$
If $k=l=\frac{q}{2}$, then $\psi$ is an eigenmap with eigenvalue $\mu=2q$.\\
\indent
Choosing $\phi=e^{i\theta}$, then it is an eigenmap from $S^1$ to $S^1$ with eigenvalue $1$. Let $g:\;{\Bbb R}^2\times{\Bbb R}^{n+1}\rightarrow {\Bbb R}^{m+1}$ be an orthogonal multiplication. Then the composition map $g(\phi,\psi):\;S^1\times S^{q-1}\rightarrow S^m$ is a bi-eigenmap with bi-eigenvalue $(1, \mu)$. Then the Hopf construction $[F_{g(\phi,\psi)}]$ is an nontrivial element in $\pi_{q+1}(S^{m+1})$. By Theorem 1.1, we have\\
 \\
{\bf Corollary 1.2}$\;\;$If there exist orthogonal multiplications $f$ and $g$, then $\forall q>2$, the homotopy class $[F_{g(\phi,\psi)}]\in\pi_{q+1}(S^{m+1})$ has a harmonic representative.\\
 \\
\indent Now the complex multiplication gives the Hopf's fibration $\psi_1:\;S^3\to S^2$, with the above $q=4,\mu=8$. The quaternion multiplication gives the eigenmap $\psi_2:\;S^7\to S^4$ with $q=8$ and $\mu=16$. Also from octonion multiplication we get the eigenmap $\psi_3:\;S^{15}\to S^8$ with $q=16$ and $\mu=32$. It is known (see [ER]) that there exist orthogonal multiplications $g_1:\; {\Bbb R}^2\times{\Bbb R}^3\to{\Bbb R}^4,\;g_2:\;{\Bbb R}^2\times{\Bbb R}^5\to{\Bbb R}^6$ and $g_3:\;{\Bbb R}^2\times{\Bbb R}^9\to{\Bbb R}^{10}$. Therefore, according to Corollary 1.2 $[F_{g_1(\phi,\psi_1)}]\in\pi_5(S^4),\;[F_{g_2(\phi,\psi_2)}]\in \pi_9(S^6)$ and $[F_{g_3(\phi,\psi_3)}]\in\pi_{17}(S^{10})$ are non trivial classes and have harmonic representatives respectively.\\
 \\
\indent
In the following sections, we will focus on the proof of Theorm 1.1.\\
\indent
 Recall that ([D]) for each $s\in(0,\frac{\pi}{2})$. There exists
 a unique $\beta_s$  which is the minimizer  of the functional
 $$J_s(\alpha) =\int_0^s(\dot{\alpha}^2+Q\cdot \sin^2\alpha)f dt$$

 over the Hilbert space
 $$X_s=\{\alpha \in H_{\rm loc}^{1}(0,s): \int_0^s(\dot{\alpha}^2
 + \alpha^2 )f dt < \infty~{\rm and}~ \alpha(s)= \frac{\pi}{2}\}.$$
 Here $$Q(t)= \frac{\lambda}{\sin^2t}+\frac{\mu}{\cos^2t} ,~~
 f(t)=\sin t\cdot \cos^qt.$$
Since $\b1_s$ is the minimizer of $J_s$, it satisfies (1.1) in
$(0,s)$. Moreover, $\b1_s(t)\to 0$ as $t\to 0$ if and only if
$J_s(\b1_s)<J_s(\f1{\pi}{2})$. In the present case, we always have
$J_s(\f1{\pi}{2})=+\infty$, it follows that \be\label{(3)}
\lim_{t\to 0}\b1_s(t)=0.\end{equation}

Similarly, we may define
 $$J_s^*(\alpha) =\int_s^{\f1{\pi}{2}}
 (\dot{\alpha}^2+Q\cdot \sin^2\alpha)f dt$$

 over the Hilbert space
 $$X_s^*=\{\alpha \in H_{\rm loc}^{1}(s,\f1{\pi}{2}): \int_s^{\f1{\pi}{2}}
 (\dot{\alpha}^2
 + \alpha^2 )f dt < \infty~{\rm and}~ \alpha(s)= \frac{\pi}{2}\}.$$
>From [D], we know that there exists a unique $\b1_s^*\in X_s^*$
which is the minimizer of $J_s^*$, satisfies (1.1) in
$(s,\f1{\pi}{2})$ and $\b1_s^*$ satisfies \be\label{(3*)}
\lim_{t\to\f1{\pi}{2}}\b1_s^*(t)=\pi \end{equation} if and only if
$$
J_s^*(\b1_s^*)=\inf_{X_s^*}J_s^*:=c_s^*<J_s^*(\f1{\pi}{2}).
$$

Since one can easily show that $c_s^*$ is uniformly bounded (by
using test functions) for $s\in (0,\f1{\pi}{2})$, while
$J_s^*(\f1{\pi}{2})\to +\infty$ as $s\to 0$, we see that
\ref{(3*)} holds true for $s>0$ small. Now we define
$$
 \alpha_s(t)=\left\{\begin{array}{clcr} \b1_s(t)&{\rm if}~~t\in (0,s]\\
\b1_s^*(t)&{\rm if}~~t\in (s,\f1{\pi}{2}). \end{array}\right.
$$
It is known ([D]) that $\al_s~:~(0,\f1{\pi}{2})\to X$ is a
continuous curve in the space
$$X=\{\alpha \in H_{\rm loc}^{1}(0,\f1{\pi}{2}): \int_0^{\f1{\pi}{2}}
(\dot{\alpha}^2
 + \alpha^2 )f dt < \infty\}.$$
 Moreover, there exists a constant $C>0$ such that
\be\label{(4)}
 J_s(\al_s)\leq C~~{\rm for}~~s\in (0,\f1{\pi}{2}).
\end{equation}
Note that $\al_s(t)$ satisfies (\ref{(1)}) in $(0,s)\cup
(s,\f1{\pi}{2})$, hence it is smooth there. However, the first
derivative $\dot{\al}_s(t)$ may have a jump at $t=s$, i.e. in
general
$$
\dot{\al}_s(s-0)\not=\dot{\al}_s(s+0).
$$
By the fact that $\al_s$ is continuous in $X$ and it satisfies
(\ref{(1)}) in $(0,s)\cup (s,\f1{\pi}{2})$, one can show that both
$\dot{\al}_s(s-0)$ and $\dot{\al}_s(s+0)$ are continuous in $
(0,\f1{\pi}{2})$. So we may define a continuous function
$$
l(s):=\dot{\al}_s(s+0)-\dot{\al}_s(s-0).
$$
It is clear that $\al_s$ is a solution to (\ref{(1)})-(\ref{(2)})
if and only if $l(s)=0$.

We notice that the arguments in [G] imply $l(s)<0$ for $s\in
(0,\f1{\pi}{2})$ with $\f1{\pi}{2}-s$ small enough. Therefore
Theorem 1.1 is a consequence of the following Proposition.
\begin{proposition}
Under the assumptions of Theorem 1.1, we have
$$
l(s)<0
$$
for $s$ sufficiently small.
\end{proposition}
The remaining part of this paper will be devoted to the proof of
the Proposition.

\section{Blow-up analysis}

We set
$$
\gamma_s(t):=\al_s(st)~~t\in(0,\f1{\pi}{s}).$$ It is clear that
$$
\ga_s(1)=\f1{\pi}{2}~~\dot{\ga}_s=s\dot{\al}_s, ~~{\rm and}
~~\ddot{\ga}_s=s^2\ddot{\al}_s.$$ By the equation (\ref{(1)}), we
have \be\label{be} \ddot{\ga}_s+s(\ctg (st)-q\tg
(st))\dot{\ga}_s-s^2(\f1{\lmd}{\sin^2(st)}+\f1{\mu}{\cos^2(st)})
\sin\ga_s\cos\ga_s=0 \end{equation} in $(0,1)\cup
(1,\f1{\pi}{s})$. It is clear that
$$
s(\ctg (st)-q\tg (st))\to \f1{1}{t},
$$
$$
s^2(\f1{\lmd}{\sin^2(st)}+\f1{\mu}{\cos^2(st)})\to\f1{\lmd}{t^2},
$$ as $s\to 0$.

Given $\e1>0$ small, there is $\d1=\d1(\e1)>0$ such that, if
$s<\d1$ the equation (\ref{be}) is well defined in
$(\e1,\e1^{-1})$ with uniformly bounded coefficients and uniformly
bounded nonlinearity. Therefore, by elliptic estimates up to
boundary (Noting that $\ga_s$ satisfies the boundary conditions
$\ga_s(\e1)\in (0, \f1{\pi}{2})$, $\ga_s(1)=\f1{\pi}{2}$, and
$\ga_s(\f1{1}{\e1})\in (\f1{\pi}{2},\pi )$.), we may get the
estimates
$$
\|\ga_s\|_{C^k([\e1 ,1])}\leq C(k,\e1),$$ and
$$
\|\ga_s\|_{C^k([1,\e1^{-1}])}\leq C(k,\e1).$$ It follows that
there exists a sequence $s_i\to 0$ such that $\ga_{s_i}\to \phi$
in $C^2([\e1,1])$ and $C^2([1,\e1^{-1}])$, where $\phi$ is $C^2$
on $[\e1,1]$ and $[1,\e1^{-1}]$, continuous on $[\e1,\e1^{-1}]$,
and satisfies the limit equation \be\label{lime}
\ddot{\phi}+\f1{1}{t}\dot{\phi}-\f1{\lmd}{t^2}\sin\phi\cos\phi
=0\end{equation} in $(\e1,1)\cup (1,\e1^{-1})$ with $
\phi(1)=\f1{\pi}{2}$.

Using a diagonal subsequence argument, we may deduce the existence
of a sequence $s_i\to 0$ such that $\ga_{s_i}\to\phi$ in
$C^2([\e1,1])$ and $C^2([1,\e1^{-1}])$ for any $\e1>0$, where
$\phi$ satisfies (\ref{lime}) in $(0,1)\cup (1,\infty)$ with $
\phi(1)=\f1{\pi}{2}$.

Consider the function $\phi$ on $[0,1]$, we have
$$
C\geq\int_0^sfQ\sin^2\al_s dt=\int_0^1sf(st)Q(st)\sin^2\ga_s dt.$$
Since the integrand of the last integral converges point wise to
$\f1{\lmd}{t}\sin^2\phi$, by Fatou's lemma we have
$$
\int_0^1\f1{\lmd}{t}\sin^2\phi dt\leq C.
$$
It follows that $\phi\not\equiv \f1{\pi}{2}$. The analysis in [D]
then shows that $\phi$ is the unique solution of the problem
$$
\left\{\begin{array}{clcr}
\ddot{\phi}+\f1{1}{t}\dot{\phi}-\f1{\lmd}{t^2}\sin\phi\cos\phi
=0\\
\phi(0+)=0,~~\phi(1)=\f1{\pi}{2}. \end{array}\right.
$$
This problem is explicitly solvable, and we get \be\label{phi}
\phi(t)=\arccos (\f1{1-t^a}{1+t^a})\end{equation} where
$a=2\sqrt{\lmd}$. One can similarly show that on $[1,\infty)$,
(\ref{phi}) is also the right expression for $\phi$.

It is useful to note that $\phi$ is a special solution of the
problem \be\label{(7)} \left\{\begin{array}{clcr}
\ddot{\phi}+\f1{1}{t}\dot{\phi}-\f1{\lmd}{t^2}\sin\phi\cos\phi
=0~~t\in (0,\infty )\\
\phi(0+)=0,~~\phi(\infty )=\pi . \end{array}\right.
\end{equation}
All solutions of (\ref{(7)}) are given by
$$
\phi_s(t)=\arccos (\f1{s^a-t^a}{s^a+t^a}),~~s\in (0,\infty ).
$$

\section{Comparison}
For comparison purposes we need to introduce a family of functions
which are solutions to the equation \be\label{(8)}
\left\{\begin{array}{clcr}&\ddot{\psi}+\ctg
t\dot{\psi}-\f1{\lmd}{\sin^2t\cos^2t}\sin\psi\cos\psi=0~~{\rm
in}~~(0,\f1{\pi}{2})\\
&\psi (0)=0,~~\psi(\f1{\pi}{2})=\pi.\end{array}\right.
\end{equation}
The solution can be explicitly given by \be \label{(10)}
\psi_s(t)=2\arctan(\ctg^{\f1{a}{2}}s\cdot\tg^{\f1{a}{2}t}),~~s\in
(0,\f1{\pi}{2})\end{equation} where, as before, $a=2\sqrt{\lmd}$.
\begin{lemma}
Let $s\in (0,\f1{\pi}{2})$ and $\psi$ be the solution of
(\ref{(8)}). Assume that for some $t_0>s$, $\al_s(t_0)>\psi
(t_0)>\max\{\t1,\f1{3\pi}{4}\}$. Then $\al_s(t)\geq\psi (t)$, for
all $t\in (t_0,\f1{\pi}{2})$. Here $$\t1=\arccos (\f1{-\lmd
(q-1)}{\mu -\lmd})\in (\f1{\pi}{2},\pi ).$$
\end{lemma}
{\it Proof.} It will be convenient to write the equation
(\ref{(1)}) as
\be\label{(11)} (f\dot{\al})\dot{}-fQ\sin\al\cos\al
=0. \end{equation}

We first show that, under the assumptions of the lemma, $\psi$
satisfies \be\label{(12)} (f\dot{\psi})\dot{}-fQ\sin\psi\cos\psi
>0~~{\rm in}~~(t_0,\f1{\pi}{2}).
\end{equation}
Note that $\psi$ as given in (\ref{(10)}) satisfies
\be\label{(13)} \dot{\psi}=\f1{\sqrt{\lmd}\sin\psi}{\sin t\cos t}.
\end{equation}
Using this we have \bs
(f\dot{\psi})\dot{}&=&\sqrt{\lmd}(\cos^{q-1}t\sin\psi )\dot{}\\
&=&\lmd\sin^{-1}t\cos^{q-2}t\sin\psi\cos\psi -\lmd (q-1)\sin
t\cos^{q-2}t\sin\psi. \es Then it is easy to derive that
$$
(f\dot{\psi})\dot{}-fQ\sin\psi\cos\psi =\sin t\cos^{q-2}t ((\lmd
-\mu )\cos\psi-\lmd (q-1))\sin\psi.
$$
>From (\ref{(13)}) we see that
$$
\psi (t)>\psi (t_0)>\t1~~{\rm for}~~t\in (t_0,\f1{\pi}{2}).
$$
Hence
$$
-\cos\psi >\f1{\lmd (q-1)}{\mu -\lmd}.
$$
It is then clear that (\ref{(12)}) holds true.

Then, we consider the function $u=\al_s-\psi$. We have
$$
u(t_0)>0~~{\rm and}~~\lim_{t\to\f1{\pi}{2}}u(t)=0.
$$
If the lemma is false, we can find $t_1\in (t_0,\f1{\pi}{2})$
where $u$ achieves a negative local minimum, i.e.
$$
u(t_1)<0,~~\dot{u}(t_1)=0,~~{\rm and}~~\ddot{u}(t_1)\geq 0.
$$
However, using (\ref{(11)}) and (\ref{(12)}), we get, at $t=t_1$,
\bs f(t_1)\ddot{u}(t_1)&<&f(t_1)Q(t_1)\left(\f1{\sin 2\al_s -\sin
2\psi}{2\al_s-2\psi}\right )u(t_1)\\ &=&f(t_1)Q(t_1)\cos\xi u(t_1)
\es where $\xi =2 (r\al_s(t_1)+(1-r)\psi(t_1))$ for some $r\in
(0,1)$. Since $\al_s(t_1),~~\psi(t_1)\in (\f1{3\pi}{4},\pi)$, we
see that $\xi\in (\f1{3\pi}{2},2\pi )$. So, $\cos\xi >0$ and we
arrive at $\ddot{u}(t_1)<0$, a contradiction! This proves the
lemma. \hfill Q.E.D.

\section{Proof of the proposition}
Consider the function $\al_s$ which satisfies (\ref{(11)}) in
$(0,s)\cup (s,\f1{\pi}{2})$. Multiplying the equation (\ref{(11)})
by $f\dot{\al}_s$ and integrating over $(0,s)$ and
$(s,\f1{\pi}{2})$ respectively we get
$$
f^2(s)\dot{\al}_s^2(s-0)=\int_0^sf^2Q\f1{\p}{\p
t}(\sin^2\al_s)dt$$
$$
-f^2(s)\dot{\al}_s^2(s+0)=\int_s^{\f1{\pi}{2}}f^2Q\f1{\p}{\p
t}(\sin^2\al_s)dt.$$ By integration by parts, we obtain
$$
f(s)(\dot{\al}_s^2(s+0)-\dot{\al}_s^2(s-0))=\int_0^{\f1{\pi}{2}}(f^2Q)'
\sin^2\al_sdt.$$ Since we [D] have $\dot{\al}_s(s+0)>0$ and
$\dot{\al}_s(s-0)>0$, we can see that
$l(s)=\dot{\al}_s(s+0)-\dot{\al}_s(s-0)>0$, if and only if the
integral on the right hand side of the last identity is positive.
Denote this integral by $I_s$, then we have \bs I_s&=&2(\mu-\lmd
q)\int_0^{\f1{\pi}{2}}\sin t\cos^{2q-1}t\sin^2\al_s dt-2\mu
(q-1)\int_0^{\f1{\pi}{2}}\sin^3 t\cos^{2q-3}t\sin^2\al_s dt\\
&:=& 2(\mu-\lmd q)I_s^1-2\mu (q-1)I_s^2. \es

To estimate $I_s^1$, we let $\ga_s(t)=\al_s(st)$ and have \bs
I_s^1&=&s\int_0^{\f1{\pi}{2s}}\sin st\cos^{2q-1}st\sin^2\ga_s dt\\
&=&s^2\int_0^{\f1{\pi}{2s}}(s^{-1}\sin st)\cos^{2q-1}st\sin^2\ga_s
dt\es Note that the integrand of the last integral converges point
wise to the function $t\sin^2\phi(t)$ as $s\to 0$, where $\phi$ is
given by (\ref{phi}). By Fatou's lemma, we have
$$
\liminf_{s\to 0}s^{-2}I_s^1\geq
\int_0^{\infty}t\sin^2\phi(t)dt=\int_0^{\infty}\f1{4t^{a+1}}{(1+t^a)^2}dt,
$$
where $a=2\sqrt{\lmd}\geq 1$. It follows that if $\lmd>1$,
\be\label{i} I_s^1\geq A(\lmd)s^2~~{\rm for}~~s~~{\rm
sufficiently}~~{\rm small}~~{\rm where}~~A(\lmd)>0,\end{equation}
and if $\lmd=1$, \be\label{ii} I_s^1\geq B(s)s^2~~{\rm
where}~~B(s)\to\infty~~{\rm as}~~s\to 0.\end{equation}

To estimate $I_s^2$, we need use the comparison lemma. First, fix
an arbitrary large number $R>0$, we have
$$
\al_s(sR)=\ga_s(R)\to\phi(R)=\arccos (\f1{1-R^a}{1+R^a})~~{\rm
as}~~s\to 0,$$ i.e. \be\label{(14)} \cos\al_s(sR)\to
-1+\f1{2}{1+R^a}.
\end{equation} Next, let $d>1$, we consider a solution of
(\ref{(8)}) as in (\ref{(10)}),
$$
\psi_{ds}(t)=2\arctan (\ctg^{\f1{a}{2}}ds\cdot\tg^{\f1{a}{2}}t).
$$
By a straightforward computation, we get
$$
\cos\psi_{ds}(Rs)=-1+\f1{2\tg^a(ds)}{\tg^a(ds)+\tg^a(Rs)}.
$$
Using the fact that
$$
\tg t=t+O(t^2)~~{\rm near}~~t=0,$$ we get \be\label{(15)}
\cos\psi_{ds}(Rs)=-1+\f1{2}{1+(R/d)^a}+O(R^2s^2).
\end{equation}
>From (\ref{(14)} and (\ref{(15)}), we can see that, if we choose
$R$ large enough, $s$ small enough,
$$
\al_s(Rs)>\psi_{ds}(Rs)\geq\pi-\e1,
$$ where $\e1>0$ can be chosen as small as we need. By the
comparison lemma, for $t\in (Rs,\f1{\pi}{2})$,
$$
\al_s(Rs)\geq\psi_{ds}(t),
$$ and consequently,
$$
\sin^2\al_s(t)\leq \sin^2\psi_{ds}(t).
$$
For simplicity, set $\e1 =\tg(ds)$. We have
$$
\sin^2\psi_{ds}=\f1{4\e1^a\tg^at}{\e1^a+\tg^a}.
$$
Then, \bs
A_s&:=&\int_{Rs}^{\f1{\pi}{2}}\sin^3t\cos^{2q-3}t\sin^2\al_sdt\\
&\leq & 4\e1^a\int_{Rs}^{\f1{\pi}{2}}\f1{\tg t}{\e1^a+\tg^at}
\sin^3t\cos^{2q-3}tdt\\
&:=&4\e1^aJ. \es Making the transformation $u=\tg t$, we get
$$
J=\int_{\tg(Rs)}^{\infty}\f1{u^{a+3}du}{(\e1^a+u^a)^2(1+u^2)^{\f1{b+5}{2}}},
$$ where $b=2q-3>-1$.
\bs J&\leq&\int_{\tg(Rs)}^1\f1{u^{a+3}du}{(\e1^a+u^a)^2}
+\int_1^\infty\f1{du}{u^{a+b+2}}\\
&:=J_1+J_2. \es Since $a+b+2>3$, we have $J_2<\infty$. Letting
$u=\e1 x$ in $J_1$, we get
$$
J_1=\int_{\tg(Rs)/\e1}^{1/\e1}\f1{\e1^{4-a}x^{a+3}dx}{(1+x^a)^2}
\leq\e1^{4-a}\int_{\tg(Rs)/\e1}^{1/\e1}x^{3-a}dx.
$$
Noting that $\e1=\tg(ds)$ and $\tg(Rs)/\e1=R/d + O(R^2s^2)$, we
see that \bs J_1&\leq &C\e1^{4-a}/R^{a-2}~~{\rm if}~~a>2,\\
J_1&\leq&C\e1^{4-a}|\log \tg(Rs)|~~{\rm if}~~a=2.\es Thus,
$$
A_s=\int_{Rs}^{\f1{\pi}{2}}\sin^3t\cos^{2q-3}t\sin^2\al_sdt \leq
C_1\e1^a+C_2\e1^{4-\d1},$$ where $\d1>0$ can be arbitrarily small.
Since $\e1=\tg(ds)$, we can rewrite the estimate as
$$
A_s=O(s^a)+O(s^{4-\d1}).
$$
Finally, we consider
$$
B_s=\int_0^{Rs}\sin^3t\cos^{2q-3}t\sin^2\al_sdt. $$
It is obvious
that, with $Rs$ sufficiently small,
$$
B_s=\tg^2(Rs)\int_0^{Rs}\sin t\cos^{2q-1}t\sin^2\al_sdt=o(1)I_s^1.
$$
Therefore, \be\label{(16)}
I_s^2=o(1)I_s^1+O(s^{2\sqrt{\lmd}})+O(s^{4-\d1}).
\end{equation}
Combining (\ref{(16)}) with (\ref{i}) and (\ref{ii}), we obtain,
$$
I_s=2(\mu-\lmd q)I_s^1-2\mu(q-1)I_s^2>0
$$
for $s>0$ sufficiently small. This proves the proposition. \hfill
Q.E.D.

\vspace{.2in}
\begin{center}
{\bf REFERENCES}
\end{center}
\vspace{.2in}

\footnotesize

\begin{description}

\item[{[D]}] W.-Y. Ding, {\em Harmonic Hopf Constructions between spheres},
Intern. J. of Math., {\bf 5} (1994), 849-860.

\item[{[ER]}] J. Eells and A. Ratto, {\it Harmonic maps and minimal immersions
with symmetries}, Annals of Math. Studies {\bf 130}, Princeton
Univ. Press, 1993.

\item[{[G]}] A. Gastel, Singularities of first kind in the harmonic map and Yang-Mills heat flows, Preprint, 2000

\end{description}

\end{document}